\documentclass[11pt]{amsart}
\usepackage{amssymb} 
\usepackage[mathscr]{eucal} 
\usepackage{xypic}   
\usepackage{amsmath} 
\usepackage[all, 2cell]{xy}
\usepackage{epsfig}
\usepackage{amscd}
\usepackage{mathrsfs}

\newtheorem{theorem}{Theorem}[section]

\newtheorem{lemma}[theorem]{Lemma}
\newtheorem{proposition}[theorem]{Proposition}

\theoremstyle{definition}
\newtheorem{definition}[theorem]{Definition}

\newtheorem{remark}[theorem]{Remark}

\numberwithin{equation}{section}

\def\XX{{\mathbb X}}
\def\YY{{\mathbb Y}}
\def\EE{{\mathbb E}}

\def\GL{\mathrm{GL}}

\def\E{{\mathfrak E}}

\begin{document}

\title[Atiyah sequences, connections and stacks]
{Atiyah sequences, connections and characteristic forms for 
principal bundles over groupoids and stacks}

\author[Indranil Biswas]
{Indranil Biswas}
\address{School of Mathematics\\
Tata Institute of Fundamental Research\\
Homi Bhabha Road, Mumbai 400005, India}
\email{indranil@@math.tifr.res.in}

\author[Frank Neumann]
{Frank Neumann}
\address{Department of Mathematics\\
University of Leicester\\
University Road, Leicester LE1 7RH, England, UK}
\email{fn8@@mcs.le.ac.uk}

\subjclass{Primary 53C08, Secondary 22A22, 58H05, 53D50}

\date{}

\keywords{groupoids, differentiable and algebraic stacks, Atiyah sequence, principal bundles, connections, curvature, de
Rham theory}

\begin{abstract}
{We construct connections and characteristic forms for principal bundles over groupoids and stacks in the differentiable, holomorphic and algebraic category using Atiyah
exact sequences associated to transversal tangential distributions.}
\end{abstract}

\maketitle
\vspace*{-0.3cm}

\section{Introduction}
\label{Intro}

In this note we outline a general theory of connections and characteristic
classes for principal $G$-bundles over groupoids and stacks, which works equally well in the differentiable, holomorphic
and algebraic setting and for groupoids, which are not necessarily \'etale.
Given a groupoid $\XX :=[X_1\rightrightarrows X_0]$ in one of these three settings
with both $X_0$ and $X_1$ smooth and the source map $s$ being a submersion, we
define a connection on it as a distribution ${\mathcal H}\,\subset\, TX_1$
transversal to the fibers of $s$. It is said to be flat if $\mathcal H$ is integrable.
Using groupoid presentations this gives rise to connections and flat
connections on stacks under certain compatibility conditions. For the particular case of
Deligne-Mumford stacks
presented by \'etale groupoids a connection always exists, because in
this case we have ${\mathcal H}= TX_1$. 

Given a Lie group $G$ and a principal $G$-bundle over a groupoid
$\XX :=[X_1\rightrightarrows X_0]$ with connection $\mathcal H$, we can define a
connection for the principal $G$-bundle. A $G$-bundle over the groupoid $\XX$ is basically
given by a principal $G$-bundle $\alpha: E_G\rightarrow X_0$ and some extra data as explained
in Section \ref{Groupoides}. A connection for the principal $G$-bundle over the groupoid $\XX$
corresponds to a splitting of the associated Atiyah sequence \cite{At}
$$0\,\rightarrow\, \text{ad}(E_G)\,\rightarrow\,\text{At}(E_G)\,\rightarrow\,
TX_0 \,\rightarrow\, 0$$
which is compatible with the data (see Section \ref{Bundles}).
This allows us to define curvature and characteristic differential forms. Using adequate
groupoid presentations this gives a general framework for studying connections and 
characteristic forms for principal $G$-bundles over stacks. In the case of Deligne-Mumford
stacks this corresponds to the constructions in \cite{BMW}.  

Full details of the results presented here will appear in \cite{BN}. There we also discuss
relations of our approach with Behrend's theory of cofoliations on stacks \cite{B1}, the various 
approaches to Chern-Weil theory for \'etale groupoids via simplicial manifolds \cite{LTX, CM, BSS} and the infinitesimal theory via Lie algebroids  \cite{CF, P}.  Another aspect of \cite{BN} is the systematic construction of secondary
characteristic classes and Deligne cohomology for principal $G$-bundles with
connections over general groupoids and stacks. This extends the theory of secondary characteristic classes of 
Cheeger-Simons differential characters to stacks and general groupoids, which for differentiable manifolds
were originally constructed in \cite{CS}, for algebraic varieties in \cite{E1} and recently in the case of \'etale groupoids in
\cite{FN}.

 \section{Groupoids, stacks and principal bundles}
 \label{Groupoides}
 
We will consider groupoids and stacks over the following three categories: the
differentiable category of ${\mathcal C}^{\infty}$-manifolds, the holomorphic category of
complex analytic manifolds and the algebraic category of smooth schemes of finite
type over a field $k$ of characteristic zero. We will refer to any of these categories
as the {\it category $\mathfrak S$ of smooth spaces}. Any smooth space $X$ has a
structure sheaf ${\mathcal O}_X$. Vector bundles are identified with their locally
free sheaves of ${\mathcal O}_X$-modules. The tangent bundle of $X$ will be denoted by $TX$.

A {\it groupoid} $\XX :=[X_1\rightrightarrows X_0]$ will mean a groupoid internal to the category $\mathfrak S$ of smooth spaces, i.e., 
the space $X_1$ of arrows and the space $X_0$ of objects 
are objects of $\mathfrak S$  and all structure morphisms
$$m: X_1\times_{X_0} X_1\rightarrow X_1, \, s,t: X_1\rightarrow
X_0,\, i: X_1\rightarrow X_1, \, e: X_0\rightarrow X_1$$ are morphisms in $\mathfrak S$ and the source map $s$ is a submersion. Here $t$ is the target map, $m$ denotes the multiplication map, $e$ the identity section and $i$ the inversion map of the groupoid.  
If in addition $s$ and $t$ are \'etale, the groupoid $\XX=[X_1\rightrightarrows X_0]$ is called {\it \'etale}.
If in addition, the {\it anchor map}
$$(s, t): X_1\rightarrow X_0\times X_0$$ is {\it proper}, the groupoid is called a {\it proper}
groupoid. 

Every groupoid $\XX=[X_1\rightrightarrows X_0]$ has an induced tangent groupoid $T\XX:=[TX_1\rightrightarrows TX_0]$. A {\it Lie group} $G$ is a groupoid $[G\rightrightarrows *]$ with one object, i.e., the space $X_0$ is just a point.

A {\it stack} $\mathfrak X$ will mean a stack over the category $\mathfrak S$ (endowed with the ``submersion'' Grothendieck topology) associated to a groupoid $[X_1\rightrightarrows X_0]$, such that $\mathfrak X$ has a {\it presentation} or {\it atlas} $X_1\rightrightarrows X_0\rightarrow {\mathfrak X}$, which is a representable surjective submersion. A stack $\mathfrak X$ is called a {\it Deligne-Mumford stack} if there is a presentation 
 $X_1\rightrightarrows X_0\rightarrow {\mathfrak X}$ such that the groupoid $[X_1\rightrightarrows X_0]$ is \'etale. Orbifolds correspond to proper Deligne-Mumford stacks (see \cite{B1, LM}). 
 Different presentations of the same stack are given by {\it Morita equivalent} groupoids
and stacks can be seen as Morita equivalence classes of groupoids \cite{BX}. Given a
groupoid $\XX$, the associated category ${\mathfrak B}\XX$ of $\XX$-torsors is a stack \cite{BX, BN}.

\begin{definition} Let $\XX=[X_1\rightrightarrows X_0]$ be a groupoid. A {\em $\XX$-space} is given by an object $P$ of $\mathfrak S$ together with a smooth map $\pi: P\rightarrow X_0$ and a map $\sigma: Q\rightarrow P, \,\,\sigma (\gamma, x):=\gamma\cdot x$, where $Q:=X_1\times_{s, X_0, \pi} P$ is the fiber product such that:
\begin{itemize}
\item[(i)] For all $x\in P$, we have: $e(\pi(x))\cdot x=x$
\item[(ii)] For all $x\in P$ and all $\gamma, \delta \in X_1$ such that
$\pi(x)=s(\gamma)$ and $t(\gamma)=s(\delta)$ we have:
$(\delta \cdot \gamma)\cdot x= \delta\cdot (\gamma \cdot x).$
\end{itemize}
Let $G$ be a Lie group and $\XX=[X_1\rightrightarrows X_0]$ a groupoid. A {\em
principal $G$-bundle} over $\XX$, denoted by $[s^*E_G\rightrightarrows E_G]$
is a principal right $G$-bundle $\pi: E_G\rightarrow
X_0$, which is also a $\XX$-space such that for all $x\in E_G$ and all
$\gamma\in X_1$ with $s(\gamma)=\pi(x)$ we have $(\gamma\cdot x)\cdot g=
\gamma\cdot (x\cdot g)$ for all $g\in G$. 
\end{definition}

It can be shown that if $\XX=[X_1\rightrightarrows X_0]$ and ${\mathbb Y}=[Y_1\rightrightarrows Y_0]$ are Morita equivalent groupoids, then there is an equivalence between their categories of principal $G$-bundles (see \cite{BX, LTX}). 

Principal $G$-bundles over a stack $\mathfrak X$ can be defined directly using a presentation. A {\em principal $G$-bundle $\E_G$ over a stack} $\mathfrak X$ is then given by a principal $G$-bundle $E_G\rightarrow X_0$ for an atlas $X_0\rightarrow {\mathfrak X}$ together with an isomorphism of the pullbacks $p_1^*E_G\stackrel{\sim}\rightarrow p_2^*E_G$ on the fiber product $X_0\times_{\mathfrak X} X_0$ satisfying the cocycle condition on $X_0\times_{\mathfrak X} X_0\times_{\mathfrak X} X_0$.
It turns out that for any submersion  $f: U\rightarrow {\mathfrak X}$ this datum defines a principal $G$-bundle $E_G\rightarrow U$ over $U$, because $X_0\times_{\mathfrak X}U\rightarrow U$ has local sections. Therefore we get a stack $\E_G$ and the $G$-multiplication map glues and there is a natural morphism $G\times {\E_G}\rightarrow {\E_G}$ (see \cite{BMW, H}). Similarly, we can define vector bundles over a groupoid $\XX$ and vector bundles over a stack $\mathfrak X$. They can be identified with respective principal $\GL_n$-bundles.

Let $G$ be a Lie group, $\mathfrak X$ a stack and $\XX=[X_1\rightrightarrows X_0]$ a groupoid presenting $\mathfrak X$. Then there is a canonical equivalence between the category of principal $G$-bundles over $\mathfrak X$ and the category of principal $G$-bundles over $\XX$ (see \cite{BX, H}). 

\section{Connections on groupoids and stacks}
\label{Connections}

Let $(\XX=[X_1\rightrightarrows X_0], s, t, m, e, i)$ be a groupoid. The source map $s\, :\, X_1\, \rightarrow\, X_0$ is a submersion and therefore the differential $ds\, :\, TX_1\, \rightarrow\, s^*TX_0$ of $s$ is surjective.  Let
$$
{\mathscr K}:=\, \text{kernel}(ds)\, \subset\, TX_1
$$
be the vertical tangent bundle for $s$. Fix a distribution given by a subbundle
${\mathscr H}\, \subset\, TX_1$
such that the natural map ${\mathscr K}\oplus {\mathscr H}\,
\rightarrow\, TX_1$ is an isomorphism. Define
$$
d_{\mathscr H}s\, :=\,(ds)\vert_{\mathscr H}\, :\, {\mathscr H}\,\rightarrow\,
s^*TX_0
$$
to be the restriction of $ds$ to ${\mathscr H}$. Note that $d_{\mathscr H}s$
is an isomorphism.

Let
$
dt\, :\, TX_1\, \rightarrow\, t^*TX_0
$
be the differential of $t$. We get a map
\begin{equation}\label{eq:theta}
\theta\, :=\,dt\circ (d_{\mathscr H}s)^{-1}\, :\,
s^*TX_0\, \rightarrow\, t^*TX_0\, .
\end{equation}
For any $y\, \in\, X_1$, let
\begin{equation}\label{eq:theta2}
\theta_y\, :=\, \theta\vert_{(s^*TX_0)_y}\, :\, T_{s(y)}X_0\,\rightarrow\,
T_{t(y)}X_0
\end{equation}
be the restriction of $\theta$ to the fiber $(s^*TX_0)_y\,=\, T_{s(y)}X_0$. We now define:

\begin{definition}\label{d-cg}
A {\em connection} on a groupoid $\XX=[X_1\rightrightarrows X_0]$ is a distribution
${\mathscr H}\, \subset\, TX_1$ such that:
\begin{itemize}
\item[(i)] for every $x\, \in\, X_0$, the image of the differential
$de(x)\,:\, T_x X_0\, \rightarrow\, T_{e(x)} X_1$ coincides with
the subspace ${\mathscr H}_{e(x)}\, \subset\,T_{e(x)} X_1$, and
\item[(ii)] for every $y\, ,z\, \in\, X_1$ with $t(y)\,=\, s(z)$,
the homomorphism $$\theta_{m(z,y)}\, :\, T_{s(y)}X_0\,\rightarrow\,
T_{t(z)}X_0$$ coincides with the composition $\theta_z\circ\theta_y$.
\end{itemize}
A connection on $\XX$ is said to be {\em flat} if the distribution $\mathscr H$ is integrable.
\end{definition}

\begin{remark}
{\rm The above definition of (flat) connections on a groupoid is equivalent to the one
given by Behrend in \cite{B1} and the one given by Laurent-Gengoux, Tu and
Xu \cite{LTX} in the case of \'etale groupoids. Flat connections in the differentiable
setting were also introduced independently as \'etalifications by Tang \cite{T}. In
fact, a connection as defined above gives a subgroupoid of the tangent groupoid
$T\XX$, which in the differentiable category is equivalent to the horizontal paths
forming a subgroupoid of the path groupoid of $\XX$. This definition is also used by
Laurent-Gengoux, Sti\'enon and Xu in \cite{LSX} to define connections for the
general framework of non-abelian differentiable gerbes via Ehresmann connections
on Lie groupoid extensions.}
\end{remark}

Let $\mathscr H$ be a connection on a groupoid $\XX=[X_1\rightrightarrows X_0]$ and
$
\tau\, :\, TX_0\, \rightarrow\, X_0
$
be the natural projection. Consider the fiber product
$
Y_1\, :=\, X_1\times_{s,X_0,\tau} TX_0\, .
$
Let
$
s'\, :\, Y_1\, \rightarrow\, TX_0
$
be the projection to the second factor and let
$
t'\, :\, Y_1\, \rightarrow\, TX_0
$
be the morphism defined by
$
t'(x\, ,v)\, =\, \theta_x(v)\,\in\, T_{t(x)}X_0 \, ,
$
where $\theta_x$ is defined as in (\ref{eq:theta2}). Furthermore, let
$
p\, :\, Y_1\, =\, X_1\times_{s,X_0,\tau} TX_0\, \rightarrow\, X_1
$
be the projection to the first factor. For
any $z\, ,y\, \in\, Y_1$ with $t'(y)\,=\, s'(z)$, define
$
m'(z, y)\,:=\, (m(p(z),p(y))\, , s'(y))\, .
$
Let
$
e'\, :\, TX_0\, \rightarrow\, Y_1
$
be the morphism defined by $v\,\mapsto\, (e(f(v))\, ,v)$ and finally let
$
i'\, :\, Y_1\, \rightarrow\, Y_1
$
be the involution defined by $(z\, ,v)\, \mapsto\, (i(z)\, ,\theta_z(v))$. We note that $([Y_1\rightrightarrows TX_0]\, , s'\, ,t'\, ,m'\, ,e'\, ,i')$
is a groupoid. In other words, $TX_0$ becomes a vector bundle over
 the groupoid $\XX$.

(Flat) connections on groupoids as defined above behave well with respect to pullbacks of cartesian morphisms and Morita equivalences of groupoids and therefore we can speak of (flat) connections on stacks. A (flat) connection on a groupoid $\XX$ 
descends to the stack ${\mathfrak B}\XX$ of $\XX$-torsors \cite{BN}.

However, even in the differentiable category these connections might not always exist. But if $\XX$ is an \'etale groupoid there is always a connection on $\XX$, because then we have ${\mathcal H}\,=\, TX_1$ and so get actually a unique flat connection (see also \cite{B1}). Therefore connections on a Deligne-Mumford stack $\mathfrak X$ always exist and can be defined using an \'etale groupoid presentation.  

\section{Connections on principal bundles over groupoids and stacks}
\label{Bundles}
Let $G$ be a Lie group and let 
$
\alpha\, :\, E_G\, \rightarrow\, X_0
$
be a principal $G$-bundle over $X_0$. The adjoint vector bundle for $E_G$ will be denoted by
$\text{ad}(E_G)$. Let
$
\text{At}(E_G)\,:=\, (\alpha_*TE_G)^G
$
be the Atiyah bundle for $E_G$ \cite{At}. We have the exact sequence of Atiyah \cite{At}:
\begin{equation}\label{eq:a-e-s}
0\,\rightarrow\, \text{ad}(E_G)\,\rightarrow\,\text{At}(E_G)\,\rightarrow\,
TX_0 \,\rightarrow\, 0\, .
\end{equation}
The projection $\text{At}(E_G)\,\rightarrow\, TX_0$ is given by
the differential $d\alpha$.
We recall that a connection on $E_G$ is a splitting of  the exact sequence (\ref{eq:a-e-s})
(see \cite{At}).

Fix a connection $\mathscr H$ on the groupoid $\XX=[X_1\rightrightarrows X_0]$ and let
$\EE_G:=([s^*E_G\rightrightarrows E_G]\, , \widetilde{s}\, ,
\widetilde{t}\, ,\widetilde{m}\, ,\widetilde{e}\, ,\widetilde{i})
$
be a principal $G$-bundle over $\XX$. We will show that
$\EE_G=[s^*E_G\rightrightarrows E_G]$ has a connection $\widetilde{\mathscr H}$ induced by $\mathscr H$. Let
$$
\phi_{X_1}\, :\, s^*E_G\, \rightarrow\, X_1
$$
be the natural projection and let
$d\phi_{X_1}\, :\, T s^*E_G\, \rightarrow\, TX_1$
be its differential. Define
$$
\widetilde{\mathscr H}\,:=\, (d\phi_{X_1})^{-1}({\mathscr H})
\, \subset\, T s^*E_G\, .
$$
It is straightforward to check that $\widetilde{\mathcal H}$ is a connection
on the groupoid $\EE_G=[s^*E_G\rightrightarrows E_G]$.

As noted above, using $\widetilde{\mathscr H}$, the tangent bundle
$TE_G$ is turned into a vector bundle over $\XX$.
We recall that $\text{At}(E_G)\,=\, (\alpha_*TE_G)^G$. Using
the action of $G$, the vector bundle $TE_G$ over
$\EE_G=[s^*E_G\rightrightarrows E_G]$ descends to make $\text{At}(E_G)$ a vector bundle
over $\XX$. All the morphisms of the vector bundle $\text{At}(E_G)$
over $\XX$ are given by descent of the corresponding morphisms of
the vector bundle $TE_G$ over $\EE_G=[s^*E_G\rightrightarrows E_G]$.

Consequently, the exact sequence (\ref{eq:a-e-s}) of vector bundles over
$X_0$ is an exact sequence of vector bundles over the groupoid $\XX$.

\begin{definition}
A {\em connection} on a principal $G$-bundle $\EE_G=[s^*E_G\rightrightarrows E_G]$
over a groupoid $\XX=[X_1\rightrightarrows X_0]$ is a splitting of the Atiyah sequence of vector
bundles over $\XX$ given by (\ref{eq:a-e-s}). It is said to be {\em flat} if it is integrable. 
\end{definition}

(Flat) connections behave again well with respect to pullbacks of cartesian morphisms and Morita equivalences of
groupoids. So we have as observed before by Laurent-Gengoux, Sti\'enon and Xu \cite{LSX} as well as by Tang \cite{T}
(see also \cite{BN}):

\begin{proposition}
If $\XX$ and $\YY$ are Morita equivalent groupoids, then there is an
equivalence between the respective categories of principal $G$-bundles with
(flat) connections over $\XX$ and $\YY$. 
\end{proposition}

Therefore we can speak of (flat) connections for principal $G$-bundles on a stack $\mathfrak X$ (see \cite{BN}).
In general these might not exist, but for principal $G$-bundles on orbifolds, i.e. proper \'etale Deligne-Mumford stacks, connections 
always exist (see also \cite{LTX}).

A stack $\mathfrak X$ determines a collection of Morita equivalent groupoids associated to the presentations of the stack and
any principal $G$-bundle ${\mathfrak E}_G$ on $\mathfrak X$ is given by principal $G$-bundles $E_G\rightarrow X_0$ over each 
atlas $X_0\rightarrow {\mathfrak X}$ satisfying appropriate glueing conditions (see for example \cite{BX}). We have the following \cite{BN}:

\begin{theorem}
Giving a (flat) connection for a principal $G$-bundle ${\mathfrak E}_G$ over a stack $\mathfrak X$ is
equivalent to giving (flat) connections for the above principal $G$-bundles on the
associated groupoids which are compatible with pullbacks along cartesian morphisms. 
Giving a (flat) connection for a principal $G$-bundle $\EE_G$ on a groupoid $\XX$ gives a (flat) connection 
for the principal $G$-bundle ${\mathfrak B}\EE_G$ of $\EE_G$-torsors over the classifying stack ${\mathfrak B}\XX$.
\end{theorem}

\section{Curvature and characteristic differential forms}
\label{Forms}

Let $\mathscr H$ be a connection on the groupoid $\XX=[X_1\rightrightarrows X_0]$. Using the
canonical decomposition of the tangent space
$$
TX_1\,=\, {\mathscr H}\oplus \text{kernel}(ds)\, ,
$$
we get a projection
$$
\wedge^j TX_1\,\rightarrow\, \wedge^j{\mathscr H}
\,\hookrightarrow\, \wedge^j TX_1\, .
$$
The composition $\wedge^j TX_1\,\rightarrow\,
\wedge^j TX_1$ gives by duality an endomorphism of the space of
$j$-forms on $X_1$.

For a differential form $\omega$ on $X_1$, the differential form on $X_1$
induced by the above endomorphism will be denoted by $H(\omega)$.

\begin{definition}
A {\em differential $j$-form} on the groupoid $[X_1\rightrightarrows X_0]$ is a differential $j$-form
$\omega$ on $X_0$ such that $H(s^*\omega)\,=\, H(t^*\omega)$.
\end{definition}

We have the following basic property:

\begin{lemma}\label{lem1}
Assume that the distribution ${\mathscr H}\, \subset\, TX_1$ is integrable.
Let $\omega$ be a differential form on $X_0$ such that
$H(s^*\omega)\,=\, H(t^*\omega)$. Then $H(s^*d\omega)\,=\, H(t^*d\omega)$.
\end{lemma}

The space of connections on a principal $G$-bundle $E_G$ over the groupoid $\XX$ is an
affine space for the space of all $\text{ad}(E_G)$-valued $1$-forms on the groupoid \cite{BN}.

Henceforth, we will assume that the distribution
${\mathscr H}\, \subset\, TX_1$ is integrable,~i.e., a {\em flat} connection.

Consider the Atiyah sequence in (\ref{eq:a-e-s}). The Lie bracket
of vector fields defines a Lie algebra structure on the sheaves
of sections of all three vector bundles.
The Lie algebra structure on the sheaf of sections of $\text{ad}(E_G)$
is linear with respect to the multiplication by functions on $X_0$, or
in other words, the fibers of $\text{ad}(E_G)$ are Lie algebras.

The Lie algebra of $G$ will be denoted by $\mathfrak g$.
Recall that $\text{ad}(E_G)=(E_G\times {\mathfrak g})/G$
is given by the adjoint action of $G$ on $\mathfrak g$. 
Since this adjoint action of $G$ preserves the Lie
algebra structure of $\mathfrak g$, it follows that the fibers
of $\text{ad}(E_G)$ are Lie algebras identified with $\mathfrak g$ up to
conjugations.

Given any splitting of  the Atiyah sequence (\ref{eq:a-e-s})
$$
D\, :\, TX_0\, \rightarrow\,\text{At}(E_G)\, ,
$$
the obstruction for $D$ to be compatible with the Lie algebra structure
is given by a section
$$
K(D)\, \in\, H^0(X_0,\, \text{ad}(E_G)\otimes\Omega^2_{X_0})\, ,
$$
which is an $\text{ad}(E_G)$-valued $2$-form on the groupoid $\XX$. The section $K(D)$
is called the {\em curvature} of the connection $D$ on the principal $G$-bundle
$\EE_G=[s^*E_G\rightrightarrows E_G]$ over the groupoid $\XX=[X_1\rightrightarrows X_0]$.
The following result allows for a general Chern-Weil theory for constructing characteristic classes \cite{BN}, which in the case of \'etale groupoids was developed before in \cite{CM} and \cite{LTX}, but where the characteristic classes are living in different cohomological spaces tailored directly for leaf spaces of foliations and \'etale groupoids.

\begin{theorem}
For any invariant form $\nu\,\in\, (\text{Sym}^k({\mathfrak g}))^G$, the closed
$2k$ form $\nu(K(D))$ on $X_0$ is a form on the groupoid.
\end{theorem}

\section*{Acknowledgements}
The second author would like to thank the Tata Institute of Fundamental Research in
Mumbai and the University of Leicester for financial support. The first author
acknowledges the support of the J. C. Bose fellowship.

\end{document}